\documentclass[s5,11pt]{article}
\usepackage{amsmath}
\newtheorem{theorem}{Theorem}

\author{ Tord Sj\"odin}
\title{Archimedian Theorems for Composite Solids}

\begin{document}

\maketitle
\begin{abstract} We consider the Center of Gravity of a solid, partly filled with some homogeneous material, and find its qualitative and quantitative properties. In particular, we prove that the Center of Gravity has its lowest position when it lies on the top surface of the material inside the solid and find a differential equation for the first moments that explains this result in both mathematical and physics terms. We make explicit calculations of this lowest position in a number of cases, such as cylinders, cones, solids of revolution, power solids, spheres and half spheres.  \end{abstract}
\paragraph{ \it AMS 2000 Subject Classsification:} Primary 26 E 25,
Secondary 46 B 20, 49 J 50  
\section{ Introduction}  The greek mathematician Archimedes was born around 287 B.C. in Syracuse on Sicily in Italy and is considered as one of the greatest scientist of all time \cite{A}, \cite{K}. See also \cite{N} and \cite{S} for more about his work. Archimedes method of exhausion enabled him to find area, volume and the Center of Gravity of many geometrical bodies. We consider a problem in this field for the Center of Gravity of a composite solid. We find the minimum point, the fixed point and the differential equation for the Center of Gravity of the partly filled solid and show that the minimum point and fixed point coincide for a large class of solids. It turns out that this problem hides a mathematical and physics principle, that can be expressed in terms of the change of rate of the first moments, that deserves to be exploited in a more general context.  
\\[0.5em]  Consider a homogeneous cylindric can with vertical axis placed on a horizontal plane $P$ partly filled with a homogeneous material. Then the Center of Gravity of the empty can is at the mid point $M$ of its symmetry axis and the same holds for the completely filled can. If we take out some of the material of a filled can then the Center of Gravity is lowered and the same happens if we start filling up an empty can. Let $h$ and $T(h)$ be the distances from the top level of the material and from the Center of Gravity of the partly filled can to the plane $P$, respectively. Then, assuming that $T(h)$ is a continuous function of $h$ on a compact interval, $T(h)$ has is a minimum point $h=h_0$. If we start filling an empty can with height $H$, then $h$ increases from $0$ to $H$, while $T(h)$ starts at $H/2$, decreases for a while and ends up at $H/2$. By continuity, there is $h=h_1$ such that $T(h_1)=h_1$ which we call a fixed point for $T(h)$. It is shown in \cite{H} that $h_0$ and $h_1$ are both unique and coincide.
\\[0.5em]
This raises two questions: a) Does the equality between the minimum point and the fixed point hold for other solids as well? b) What is the underlying theory in mathematics and/or physics that explains why minimum points and fixed points coincide? We start by answering question a) for a number of standard solids and then find that the answer to question b) lies in a differential equation for the  first moment of the partly filled solid. This result remains true for more general solids, as long as the quantities involved are well defind.\\[0.5em] The plan of this paper
is as follows. Section 2 solves the problem for the cylinder (Theorem 1), while Sections 3 and 4 treat the cone and a general solid of revolution (Theorem 2). Sections 5 and 6 consider power solids, spheres and half spheres. In section 7 we find the differential equation (Theorem 3) that governs the behaviour of the Center of Gravity and explains Theorem 2 in terms of the rate of change of the first moment of the partly filled solid (Theorem 4). We conclude with some comments and historical background.\section { The cylinder}  We start with a mathematical treatment of the cylinder that follows the analysis in \cite{H}. An open cylinder with height $H$ is placed vertically on a horizontal plane $P$and filled with a homogeneous material up to the level $h$ from $P$. Then the Center of Gravity of the cylinder and its content lies on the central axis of the cylinder at the distance $T_0(h)$ from $P$. By the homogenety of the cylinder and the material we have $T_0(0)=T_0(H)=H/2$. We want to find the lowest point of the Center of Gravity, i.e.  the minimum value of the function $T_0(h)$ over the interval $0\leq h\leq H$.
  Assume that the weight of the empty cylinder is $M$ and that the weight of the material in a completely filled cylinder is $m$. A straight forward calculation gives that 
\begin{displaymath}T_0(h)=\frac{M\cdot \frac{H}{2}+m\cdot \frac{h}{H}\cdot \frac{h}{2}}{M+m\cdot \frac{h}{H}}=\frac{1}{2}\cdot \frac{M\cdot H^2+m\cdot h^2}{M\cdot H+m\cdot h},\, 0\leq h\leq H,\end{displaymath}
with derivative
\begin{displaymath}T_0^\prime (h)=\frac{m}{2}\cdot \frac{m h^2+2MHh-MH^2}{(M H+m  h)^2},\, 0\leq h\leq H.\end{displaymath}
It follows that $T_0(h)$ has a unique minimum point in [0,H] given by the positive root of the quadratic equation
\begin{equation}m h^2+2MHh-MH^2=0.\end{equation}
We can now calculate the minimum value of $T_0(h)$, but (1) does not explicitely show the position of the Center of Gravity relative to the top surface of the material. The continuous map  $T_0:[0,H]\rightarrow [0,H]$ has a fixed point $h_1$ given by $T_0(h_1)=h_1$ and it turns out that also $h_1$ is given by (1). This common value of the minimum point and the fixed point is obtained from the formula  
\begin{displaymath}h= -\frac{MH}{m}+\sqrt{\frac{M^2H^2}{m^2}+\frac{MH^2}{m}}=\frac{MH}{m}\big( \sqrt{1+\frac{m}{M} }-1\big). \end{displaymath} 
 \begin{theorem}An open cylinder with finite height is placed vertically on a horizontal plane and is partly filled with some homogeneous material. Then the Center of Gravity of the cylinder with its content is at its lowest position if and only if it lies on the top surface of the material inside the cylinder. The value of this position is given by (1).\end{theorem}
 This nice result has another and simpler proof  \cite{H}. Turn the cylinder horizontal and such that the material inside it keeps its shape. Assume that the material inside the cylinder is at the position described in the theorem and support the cylinder from below at that point. Then, increasing or decreasing the amount of material will move the Center of Gravity away from the bottom of the cylinder and the theorem follows. We will put this problem in a wider context and consider a number of other solids, including a general solid of revolution, with a method based on slicing (Cavalieri's principle).
 \section{The cone} Consider a top down cone with vertical axis and height $H$ placed with its vertex on a horizontal plane $P$. The weights of the empty cone and the material of a completely filled cone are $M$ and $m$ respectively. Let $T_1(h)$ denote the distance from the Center of Gravity to $P$ of a cone filled with material up to the level at the distance $h$ from $P$. A straight forward calculation finds that
\begin{displaymath}T_1(h)=\frac{M\cdot \frac{3}{4}\cdot H+m\cdot (\frac{h}{H})^2\cdot \frac{2}{3}\cdot h}{M+m\cdot (\frac{h}{H})^2}=\frac{1}{12}\cdot \frac{9MH^3+8mh^3}{MH^2+mh^2},\, 0\leq h\leq H,\end{displaymath}
and
\begin{displaymath}T_1^\prime (h)=\frac{mh}{12}\cdot \frac{8mh^3+24MH^2h-18MH^3}{(MH^2+mh^2)^2},\, 0\leq h\leq H.\end{displaymath}
It is easy to see that $T^\prime_1(h)=0$ has a unique positive solution in [0,H] given by the third degree polynomial equation
\begin{equation}4mh^3+12MH^2h-9MH^3=0.\end{equation}
It should now be no surprise that the equation (2) also gives the unique fixed point of the map $T_1:[0,H]\rightarrow [0,H]$, which means that Theorem 1 holds for the cone as well.\\[0.5em] {\it Remark.} One difference between  the cylinder and the cone is that the empty cylinder and the material of a filled cylinder have the same Center of Gravity, which is not the case for the cone. The Center of Gravity of the empty cone has distance $3H/4$ from the vertex and the Center of Gravity of the material in the filled cone lies $2H/3$ from the vertex. This does however not affect the fact that minimum points and fixed points coincide.
\section{Solids of revolution} The cylinder and the cone can both be generated by rotating a straight line around an axis. We will consider solids generated by other curves, such as power functions and circles. We begin with a general curve and get the various special cases from that. Define a plane set in $R^3$ by
$ D=\{ (x,0,z);0\leq x\leq g(z) \, \textrm{and}\, 0\leq z\leq H\}$, where $g(z)$ is a positive and continuously differentiable function, $0\leq z\leq H$. Let $B$ denote the solid with surface $S$ that is generated when the set $D$ and the curve $x=g(z)$ are rotated around the vertical $z-$axis. Denote the plane $z=h$ by $P_h$. We use Cavalieri's principle to calculate the areas, volumes and moments we need to find the Center of Gravity of $B$. The volume $V$ of $B$ and the area $A$ of $S$ are given by
\begin{displaymath}V=\pi \int\limits _0^H g(z)^2\, dz\quad \textrm{and}\quad A=2\pi \int\limits _0^H g(z)\cdot \sqrt{1+g^\prime (z)^2}\, dz,
\end{displaymath}
respectively. Assume that the surface has area density $\alpha$ and that the material has volume density $\beta$. Then the zero--th moment $S_0$ and the first moment $S_1$ of $S$, relative to $P_0$, are given by
\begin{displaymath}S_0=2\pi  \int\limits _0^H g(z)\cdot \sqrt{1+g^\prime (z)^2}\, dz\quad \textrm{and}\quad S_1=2\pi  \int\limits _0^H z\cdot g(z)\cdot \sqrt{1+g^\prime (z)^2}\, dz.
\end{displaymath}
  This gives us the following formula for the distance $T_{sr}(h)$ between the Center of Gravity of the solid $B$ filled with material between the planes $P_0$ and $P_h$ as
\begin{equation}T_{sr}(h)=\frac{m_1(h)}{m_0(h)}=\frac{\alpha\cdot S_1+\pi\beta\cdot \int\limits _0^h z\cdot  g(z)^2\, dz}{\alpha\cdot S_0+\pi\beta\cdot \int\limits _0^h  g(z)^2\, dz},
\end{equation}
where $m_0(h)$ and $m_1(h)$ are the zero--th and first moments of the solid and its content relative to $P_0$.
Further,
\begin{displaymath}T^\prime_{sr}(h)=\frac{\pi\beta\cdot g(h)^2}{ m_0(h)^2}\cdot  \big( \alpha S_0h+\pi\beta h\cdot \int\limits _0^h   g(z)^2\, dz-\alpha S_1-\pi\beta\cdot \int\limits _0^h z\cdot g(z)^2\, dz \big)
\end{displaymath}
and \begin{displaymath}T_{sr}(h)-h=\frac{ \alpha S_0h+\pi\beta h\cdot \int\limits _0^h   g(y)^2\, dy-\alpha S_1-\pi\beta\cdot \int\limits _0^h y\cdot g(y)^2\, dy}{m_0(h)}
\end{displaymath}
Clearly, $T_{sr}:[0,H]\rightarrow [0,H]$, since $S_1\leq H\cdot S_0$, and therefore $T_{sr}$ has a fixed point. Let $F(h)$ denote the last factor in $T_{sr}^\prime (h)$, then $F(0)=-\alpha S_1<0$ and $F(H)\geq 0$. Since $F^\prime (h)>0$ we conclude that $F(h)$ has a unique zero in $[0,H]$. It follows that $T_{sr}(h)$ has a unique fixed point and $T_{sr}^\prime(h)=0$ if and only if $T_{sr}(h)=h$, where $h$ is the unique solution of the equation
\begin{equation}\alpha S_0h+\pi\beta h\cdot \int\limits _0^h   g(y)^2\, dy-\alpha S_1-\pi\beta\cdot \int\limits _0^h y\cdot g(y)^2\, dy =0\end{equation} in the interval $[0,H]$. This proves the following theorem stating that fixed points and minimum points agree for solids of revolution.
\begin{theorem} Let $B$ be the solid of revolution described in this section filled with a homogeneous materiel between the planes $P_0$ and $P_h$. Then the Center of Gravity of the solid with its content is at its lowest position if an only if it lies on the top surface of the material and its value is calculated by (4).\end{theorem}
\section{Power solids.}
We choose $g(z)$ as the power function $g(z)=z^p$, $0\leq z\leq H$, for some $p>0$, and call the corresponding solid $B$ a power solid (ps). A straight forward calculation gives the expression
\begin{displaymath}T_{ps}(h)=\frac{\alpha\cdot S_1+\pi\beta\cdot h^{2p+2}/(2p+2)}{\alpha\cdot S_0+\pi\beta\cdot h^{2p+1}/(2p+1)}\end{displaymath}
for the distance between the Center of Gravity of $B$ and the plane $P_0$. The unique minimum point and fixed point of $T_{ps}(h)$ is given by the equation
\begin{equation}\alpha S_0h-\alpha S_1+\pi\beta\cdot \frac{h^{2p+2}}{(2p+1)(2p+2)}=0,\end{equation}
where
\begin{displaymath}S_0=2\pi\cdot \int\limits _0^H\sqrt{1+p^2\cdot z^{2p-1}}\, dz \quad \textrm{and}\quad  S_1=2\pi\cdot \int\limits _0^Hz\cdot \sqrt{1+p^2\cdot z^{2p-1}}\, dz .
\end{displaymath}
\section{Spheres and half spheres}We start with a sphere of radius $R$ centered at $(0,0,R)$ and generated by the curve $x=g(z)=\sqrt{R^2-(R-z)^2}$, $0\leq z\leq 2R$. We find that $g(z)\cdot \sqrt{1+g^\prime (z)^2}=R$ and hence $S_0=4\pi R^2$ and $S_1=4\pi R^3$. This gives the following formula for the distance $T_{sp}(h)$ between the Center of Gravity of the partly filled sphere and the plane $P_0$
\begin{displaymath}T_{sp}(h)=\frac{4\alpha\cdot R^3+\beta\cdot (2Rh^3/3-h^4/4)}{4\alpha\cdot R^2+\beta\cdot (Rh^2-h^3/3)},\quad 0\leq h\leq 2R,
\end{displaymath}
and the unique minimum/fixed point is given by the fourth order polynomial equation
\begin{equation}4\alpha\cdot R^3+\beta\cdot (2Rh^3/3-h^4/4)-4\alpha\cdot R^2h-\beta\cdot (Rh^2-h^3/3)h=0.\end{equation}
A simple check shows that $T_{sp}(0)=T_{sp}(2R)=R$, as expected, since  both the material in a completely filled sphere and the sphere itself have Center of Gravity at the center of the sphere. Solving (6.1) for $\alpha$ we get $\alpha =(\beta /48R^2)\cdot$ $h^3(4R-h)/(R-h)$, $0\leq h<R$, which for fixed $R$ and $\beta$ defines $\alpha$ as a continuous, strictly increasing and convex function of $h$, Taking inverse, $h$ is a strictly increasing and concave function of $\alpha$, for $\alpha >0$,
 \\[0.5em] The half sphere with radius $R$ and centered at $(0,0,R)$ is generated by the curve $x=g(z)=\sqrt{R^2-(R-z)^2}$, $0\leq z\leq R$, and has $S_0=2\pi R^2$ and $S_1=\pi R^3$. The distance $T_{hsp}(h)$ between the Center of Gravity of the partly filled half sphere and the plane $P_0$ is then
\begin{displaymath}T_{hsp}(h)=\frac{ \alpha\cdot R^3+\beta\cdot (2Rh^3/3-h^4/4)}{2\alpha\cdot R^2+\beta\cdot (Rh^2-h^3/3)},\quad 0\leq h\leq R,\end{displaymath}
and the unique minimum/fixed point satisfies
\begin{equation}  \alpha\cdot R^3+\beta\cdot (2Rh^3/3-h^4/4)-2\alpha\cdot R^2\cdot h-\beta\cdot (Rh^2-h^3/3)\cdot h=0.\end{equation}
Here we have $T_{hsp}(0)=R/2$, while $T_{hsp}(R)$ depends on $R$, $\alpha$ and $\beta$, since the Centers of Gravity of the half sphere is not the same as the Center of Gravity of the material in a completely filled half sphere. Solving (6.2) for $\alpha$ gives $\alpha =(\beta  /12R^2)\cdot (h^4-4Rh^3)/(2h-R)$, $0\leq h<R/2$, which defines $\alpha$ as continuous and strictly increasing function of $h$. Taking inverse, $h$ is a strictly increasing function of  $\alpha$, for $\alpha\geq 0.$ 
\section{Differential equations}
The two proofs of Theorem 1 gives an analytic and a geometric explanation why minimum points and fixed points coincide using first moments. We will give another explanation in terms of a differential equation. Our calculations above show that the functions $T_i (h)$, $i=0,1$ and $T_{sr}(h)$, in Sections 2, 3 and 4 satisfy the first order ordinary linear differential equations
\begin{equation}T_0^\prime (h)=\frac {-4}{MH+mh}\cdot (T_0(h)-h),\end{equation}
 \begin{equation}T_1^\prime (h)=\frac {-m}{3}\cdot \frac{h}{MH^2+mh^2}\cdot (T_1(h)-h),\end{equation}
 and \begin{equation}T_{sr}^\prime (h)=\frac {-\pi\beta\cdot g(h)^2}{\alpha\cdot S_0+\pi\beta\cdot \int\limits _0^h  g(z)^2\, dz}\cdot (T_{sr}(h)-h),
 \end{equation}
respectively. The differential equation (10) explains why minimum points and fixed points coincide for all solids of revolution, as stated in Theorem 2. All the  equations (8), (9) and (10) are of the form $T^\prime (h)=-\phi (h)\cdot (T(h)-h)$, for some positive function $\phi (h)$. We want to interpret the function $\phi (h)$ in terms of physical quantities of the solid $B$. For this purpose, we rewrite (10) as
\begin{equation}T_{sr}^\prime (h)=-\frac {m_0^\prime (h)}{m_0(h)}\cdot (T_{sr}(h)-h),\quad  m_0(h)=\alpha\cdot S_0+\pi\beta\cdot \int\limits _0^h  g(z)^2\, dz,\end{equation}
where $m_0(h)$ is the zero-th moment of the partly filled solid.
  Equation (11) expresses the increase/decrease of $T_{sr}(h)$ as a product of the distance between the Center of Gravity $T_{sr}(h)$ and the top surface of the material inside the solid $B$ and the relative change of $m_0(h)$.  
   \begin{theorem} Let $B$ be the solid of revolution in Section 4 and let $T_{sr}(h)$ be the distance between $P_0$ and the Center of Gravity of $B$ filled with material upp to the plane $P_h$. Then $T_{sr}(h)$ is defined by (3), the minimum point of $T_{sr}(h)$ coincides with the fixed point and $T_{sr}(h)$ satisfies the differential equation (11).
\end{theorem}
Equation (11) explains why minimum points and fixed points coincides in mathematical terms. We also would like to have an explanation in physical quantities of the solid. To achieve this we rewrite (11) as
\begin{equation}\big( T_{sr}(h)\cdot m_0(h)\big) ^\prime =m_0^\prime (h)\cdot h.\end{equation}
Here, $ T_{sr}(h)\cdot m_0(h)$ is the first moment of the partly filled solid and the right hand side of (12) is the rate of change of this moment. Thus (12) is nothing more than the rate of change of the first moment expressed in two different ways. This is the principle alluded to in the introduction and is our final explanation why Theorems 2 and 3 hold.\begin{theorem} Under the assumptions of Theorems 2 and 3, the equality of the minimum point and the fixed point can be reduced to the identity (12) for the rate of change of the first moment of the partially filled solid.\end{theorem}
Theorems 3 and 4 hold for more general solids as well. Let $G$ be an open set in $R^3$ with sufficiently smooth boundary $S$ lying between two horizontal planes $P_0$ and $P_H$ and put $B=G\bigcup S$. Assume that $S$ has well defined zero--th moment $S_0$ and first moment $S_1$, relative to $P_0$. Let $f(h)$ be the area of the set $B\bigcap  P_h$ and assume that $f(h)$ is continuous for $0\leq h\leq H$. Then the Center of Gravity of $B$ filled with material between $P_0$ and $P_h$ is given by
\begin{displaymath}T(h)=\frac{\alpha\cdot S_1+ \beta\cdot \int\limits _0^h z\cdot  f(z)\, dz}{\alpha\cdot S_0+ \beta\cdot \int\limits _0^h  f(z)\, dz},\quad 0\leq h\leq H,
\end{displaymath}
and one readily verifies that $T^\prime (h)=-m^\prime (h)/m(h)\cdot (T(h)-h)$, where $m(h)=\alpha\cdot S_0+\pi\beta\cdot \int\limits _0^h  f(z)\, dz$ and $\alpha ,\beta$ are as in Section 4. If we let $S_1(h)=T(h)\cdot m(h)$ denote the first moment of the partly filled solid relative to the plane $P_0$ then the equation becomes $S_1^\prime (h)=m^\prime (h)\cdot h$, which expresses the change of the first moment in two ways. The explicit position of this point is given by the equation $T(h)=h$ and depends on the geometry of the solid,  the area density $\alpha$ of the surface and the volume density $\beta$ of the material.
\section{Concluding remarks} Archimedes was  
the first mathematician to define and study equilibrium of masses and Centers of Gravity, see \cite{A} and \cite{K} for a detailed history and \cite{N} for a translation of some of his work. He mostly dealt with homogeneous bodies in the plane and in space and it is not known to us if he met the problem described in this paper. The  second proof of Theorem 1 could very well have been known to Archimedes and his contemporaries and possibly been lost somewhere in the history.

Address: Tord Sj\"odin, Department of Mathematics and Mathematical Statistics, University of 
Ume\aa , 
 S-901 87  Ume\aa , Sweden. E-mail: tord.sjodin@math. umu.se

\end{document}